\documentstyle[epsf]{article}

\newcommand{\ee}{\end{equation}}
\newcommand{\be}{\begin{equation}}

\begin{document}

\centerline{ \huge \bf An Asymptotic Reduction of}
\vskip 0.5 cm 
\centerline{ \huge \bf a Painlev\'e VI equation to a  Painlev\'e III}

\vskip 0.5 cm 
\centerline{\huge Davide Guzzetti}
~~~~~~~~~~~~~~~~~~~~~~~~~~~~~~~~~~~~~~~~~~~~~~~~~~~~~~~~~~~~~~~~~~~~~~~~~~~~~~~~~~~\footnote{ 1) International School of Advanced Study SISSA/ISAS, Trieste, Italy.

~~2)  Korea Institute of Advanced Study KIAS, Seoul, South Korea.   
 
 ~~~~~~~E-mail: davide$\underline{~~}$guzzetti@yahoo.com}
\vskip 1 cm
\noindent
{\bf Abstract:} When the independent variable is close to a
critical point, it is  shown that PVI can be asymptotically reduced to
PIII. In this way, it is possible  to compute the leading term of the critical behaviors of PVI
transcendents starting from the behaviors of PIII transcendents.

\vskip 1 cm

\section{Introduction} 
As it is well known, the first five Painlev\'e equations PI, PII, PII, PIV, PV
 can be
obtained from the sixth PVI, by a step-by-step degeneration process
\cite{Ince}. 

Here we present a different reduction of PVI to PIII. When the
 independent variable is close to a critical point, we  show that it is possible to reduce  
 PVI, with $\alpha=(2\mu-1)^2/2\in{\bf C}$, $\beta=\gamma=0$,
 $\delta=1/2$, namely:
$$
y_{ss} ={1\over 2}\left[ 
{1\over y}+{1\over y-1}+{1\over y-s}
\right]
           \left(y_s\right)^2
-\left[
{1\over s}+{1\over s-1}+{1\over y-s}
\right]{y_s}
+~~~~~~~~~~~~~~~~~~~~~~~~~~~~~~~~~~~~
$$
\be
\label{p6}
~~~~~~~~~~~~~~~~~~~~~~~~~~~~~~~~~~~~~~~~~~~~~~
+{y(y-1)(y-s)\over 2s^2 (s-1)^2}
\left[
(2\mu-1)^2+{s(s-1)\over (y-s)^2}
\right], 
\ee
 to the following PIII, with $\alpha=\beta=0$, $\gamma=-\delta=1$, namely:
\be
\label{p3}
\hat{y}_{\theta\theta}={1\over y} (\hat{y}_\theta)^2-{1\over \theta}
\hat{y}_\theta + \hat{y}^3-{1\over \hat{y}}.
\ee
 The notation $y_{*}$ stands for $dy/d*$. 
  We consider here, for definiteness, the case of the critical point $s=0$. The reduction will be done accordingly for $s\to 0$ (and
 $s\sim \theta^2$). The convergence is intended for bounded arg$(s)$.  
  This is an {\it asymptotic
reduction}. We show that, remarkably, it allows to  reproduce the
correct leading term of the critical behavior of PVI 
transcendents, which are  classified in \cite{guz2010}. 
 The equation ({\ref{p6}) is important in the theory of semi-simple
 Frobenius manifolds of dimension 3 (see \cite{Dub1} \cite{Dub2}).

\subsection{A $3\times 3$ isomonodromy representation} 
 Usually, PVI is regarded as the isomonodromy deformation equation for a $2\times 2$ fuchsian system
 with four singularities \cite{SMW}. Here we regard (\ref{p6})  as the
 isomonodromy deformation condition of the $3\times 3$ linear system
 (L1) below, having  a fuchsian singularity at $z=0$ and an irregular singularity of
 rank one at $z=\infty$. This system for the general PVI is described in \cite{Marta}. For the particular PVI we consider here, related to Frobenius manifolds, (L1) is introduced and studied in \cite{Dub1} \cite{Dub2},
 where  the following $3\times 3$ Lax Pair is given: 
$$
 {\rm (L1)}:~~{dY\over dz}= \left[U+{V\over z}\right]Y,~~~~~{\rm (L2)}:~~{\partial Y\over \partial u_i}= 
\bigl[zE_i+V_i\bigr]Y,~~~~~i=1,2,3.
$$ 
The $3\times 3$ matrix coefficients are:
$$
 U=\hbox{diag}(u_1,u_2,u_3).
$$
$$
V=V(u_1,u_2,u_3),~~~V^T=-V, ~~~~~V\hbox{
   has diagonal form = diag}(\mu,0,-\mu).
$$
$$
(E_k)_{kk}=1,~~~~~ (E_k)_{ij}=0, ~~i,j\neq k.~~~~~
(V_k)_{ij}={\delta_{ki}-\delta_{kj}\over u_i-u_j} V_{ij}.
$$
A fundamental solution of (L1) at $z=0$ has representation  $Y(z)=
\left[\sum_{p=0}^\infty \phi_p(u)
  z^p\right]~z^{\hbox{diag}(\mu,0,-\mu)}
 z^\Lambda$, where
$\Lambda_{ij}=0$ if $\mu_i-\mu_j\neq n>0$, $n\in {\bf Z}$. 
The necessary and sufficient
condition for the dependence of the system (L1) on $(u_1,u_2,u_3)$ to
be isomonodromic is \cite{SMW}: 
$$
 {\partial V\over \partial u_i}= [V_i,V],~~~{\partial \phi_0\over \partial u_i}= V_i
 \phi_0.
$$
The equations $\sum_i {\partial V \over \partial u_i}=0$ and $\sum_i
u_i  {\partial V \over \partial u_i}=0$, imply  that $
V=V(s)$, where  $s={u_3-u_1\over u_2-u_1}$.  
Thus, if we let: 
 \be
   V(s)= \pmatrix{  0  &   -\Omega_3     &    \Omega_2  \cr
                    \Omega_3  &  0   &  -\Omega_1  \cr
                   -\Omega_2  & \Omega_1  &  0   \cr
                  }
, 
\label{V}
\ee
the equation for $V$ becomes: 
\be
\left\{ \matrix{
{d \Omega_1\over ds} = {1 \over s} ~ \Omega_2 \Omega_3\cr\cr
{d \Omega_2\over ds} = {1 \over 1-s} ~ \Omega_1 \Omega_3 \cr\cr
{d \Omega_3\over ds} = {1 \over s( s-1)} ~ \Omega_1 \Omega_2  
}
\right.
\label{system}
\ee

$$\Omega_1^2+\Omega_2^2+\Omega_3^2=-\mu^2$$
 This system is equivalent to (\ref{p6}), and the following holds [see  \cite{guz2001} for details]: 
\be 
  y(s)= { - s A(s) \over 1 - s (1 +A(s))},~~~~ A(s) := \left[ {\Omega_1
  \Omega_2 + \mu \Omega_3 \over \mu^2 + 
  \Omega_2^2} \right]^2 = \left[ {\Omega_1
  \Omega_2 + \mu \Omega_3 \over \Omega_1^2+\Omega_3^2} \right]^2
\label{PVI}
\ee
This last equation allows to compute the critical behavior of $y(s)$
if we know that of the $\Omega_j(s)$'s.


\vskip 0.3 cm 
\section{Asymptotic Reduction of PVI to PIII}
 
 As $s\to 0$ we do the following approximation to the system
 (\ref{system}): 
$$
{d \Omega_1^{(a)}\over ds} = {1 \over s} ~ \Omega_2^{(a)} \Omega_3^{(a)},~~~
{d \Omega_2^{(a)}\over ds} =  ~ \Omega_1^{(a)} \Omega_3^{(a)},~~~
{d \Omega_3^{(a)}\over ds} = -{1 \over s} ~ \Omega_1^{(a)} \Omega_2^{(a)} . 
$$
with 
 $\Omega_j(s)\sim \Omega_j^{(a)}(s)$ for  $s\to 0$. The superscript $(a)$ stands for "asymptotic". The reduced system
has a 
first integral:  
$$ 
  (\Omega_1^{(a)})^2+(\Omega_3^{(a)})^2 = R^2 \in{\bf C}
$$
This implies that we can introduce the new dependent variable
$\phi(s)$ as follows:
$$
   \Omega_3^{(a)}=R \sin \phi,~~~~\Omega_1^{(a)}=R \cos \phi,~~~R\neq 0
$$
The system becomes:
$$
{d \phi \over d s} = -{1\over s} \Omega_2^{(a)},~~~
{d \Omega_2^{(a)} \over d s} = {R^2 \over 2} \sin(2 \phi).
$$
Thus:
$$
 {d^2\phi\over ds^2} +{1\over s} {d\phi\over ds} +{R^2\over 2s}\sin
 \phi=0
$$
With another change of variables:  
 $$ 
   u := 2 \phi, ~~~~~ s= {x^2\over 4 R^2}, ~~~(\hbox{namely } 2\phi(s)=u(2R\sqrt{s}))
$$ 
 we obtain the following particular form of the Painlev\'e III
equation:
\be
    u_{xx} +{1\over x} u_x + \sin(u)=0.
\label{PIII}
\ee
A last change of variables is necessary: 
$$
    x=2i\theta,~~~~\hat{y}= \exp\left({i u\over 2} \right).  
$$
This gives the PIII equation in standard form (\ref{p3}): 
$$
  \hat{y}_{\theta\theta}={1\over \hat{y}} (\hat{y}_{\theta})^2-{1\over \theta}
  \hat{y}_{\theta}+\hat{y}^3-{1\over 
  \hat{y}} 
$$
To summarize, the change of variables is:
$$
 \Omega_1^{(a)}={R\over 2}\Bigl(\hat{y}^{-1}+\hat{y} \Bigr),~~~ 
\Omega_3^{(a)}={iR\over 2}\Bigl(\hat{y}^{-1}-\hat{y} \Bigr),~~~
\Omega_2^{(a)}= i s {d\hat{y}\over ds}
$$
$$
   \hat{y}(s):=
   \hat{y}\left(\theta(s)\right)=\hat{y}\left({R\sqrt{s}\over i}\right).
$$

\section{From asymptotic behaviors of PIII to behaviors of PVI}

Let $s\to 0$, $|\arg x |<\pi$. In   \cite{guz2010} we
classified the critical behaviors of the PVI transcendents into a few classes, in the case when there is a one to one correspondence between branches of PVI-transcendents and points in the space of the associated monodromy data.  
The critical behaviors  are decided by the value of a complex 
``exponent'' $\sigma$ such that $0\leq \Re \sigma \leq 1$. Let $\nu$
be a real number. Let also $a\neq 0$ and $C$ be two complex numbers, which, 
together with $\sigma$, play the role of constants of integration. According to \cite{guz2010},  the equation (\ref{p6}) has solutions with branches admitting the the following critical behaviors for $s\to 0$:

\vskip 0.2 cm 
\noindent
{\bf 1) Small-power-type behaviors (Jimbo \cite{Jimbo}) -- 4 real parameters:}
$$
 y(s)= a x^{1-\sigma}(1+O(s^{\sigma}+s^{1-\sigma})) , ~~~0<\Re\sigma<1.~
$$
 {\bf 2) Sine-type oscillatory behaviors -- 3 real parameters:}
$$
 y(s)= s \left[\sin^2\left(\nu\ln s + C\right) +O(s) \right]
,~~~\sigma =2i\nu.
$$
{\bf 3) Inverse sine-type oscillatory behaviors -- 3 real parameters:}
$$
y(s)={1\over 1-{4\nu^2+(2\mu -1)^2\over 4 \nu^2} \sin^2(\nu \ln s +C)
  +O(s)},~~~\sigma=1+2i\nu.
$$
{\bf 4) Log-type behaviors and Taylor expansions-- 2 real parameters:}
$$
y(s)= -{4\over (2\mu-1)^2} {1\over (\ln s +C)^2}\left[1+O\left({1\over
  \ln^2 x}  \right)\right],~~~\sigma=1.
$$
$$
y(x)=as+O(s^2),~~~\sigma=0.~~~~~~~~~~~~~~~~~~~~~~~~~~~~~~~~~~~~~
$$
The last is a convergent Taylor series, degeneration of a log-behavior which
occurs  when
$\beta={1\over \delta}-1=0$ in the general PVI. The higher order terms
in {\bf 1), 2), 3)} can be written as convergent expansions in $s$ and
$s^\sigma$, as it is explained in \cite{guz2010}. 

\vskip 0.3 cm 
We are going to show  that we can obtain the leading term of the above
behaviors, when we substitute into the $\Omega_j^{(a)}$'s the asymptotic
expansions of the solutions of (\ref{p3}), according to the formulas:

\be
\label{bikkuri}  
  y(s)= { - s A(s) \over 1 - s (1 +A(s))},~~~~ A(s) \sim \left[ {\Omega_1^{(a)}
  \Omega_2^{(a)} + \mu \Omega_3^{(a)} \over \mu^2 + 
  (\Omega_2^{(a)})^2} \right]^2 = \left[ {\Omega_1^{(a)}
  \Omega_2^{(a)} + \mu \Omega_3^{(a)} \over
  (\Omega_1^{(a)})^2+(\Omega_3^{(a)})^2}
 \right]^2,~~~s\to 0. 
\ee
and:
$$
 \Omega_1^{(a)}={R\over 2}\Bigl(\hat{y}^{-1}+\hat{y} \Bigr),~~~ 
\Omega_3^{(a)}={iR\over 2}\Bigl(\hat{y}^{-1}-\hat{y} \Bigr),~~~
\Omega_2^{(a)}= i s {d\hat{y}\over ds}.
$$

\vskip 0.2 cm 
 At this point, we need the asymptotic behaviors of the solutions of the PIII equation  (\ref{p3}). We can find them in  
 \cite{Tracy}.  

\vskip 0.2 cm 
{\bf i)} The first expansion  we consider is (already in variable $s$):
$$ 
  \hat{y}(s) =B s^{\sigma\over 2} \left( 
              1 + {1\over 4} \left({R\over B} \right)^2{s^{1-\sigma} \over
              (1-\sigma)^2} -{1\over 4} {(RB)^2\over (1+\sigma)^2} s^{1+\sigma} +
\sum_{j=3}^{+\infty}\sum_{k=2}^{j+1} c_{jk} s^{{1\over 2}( j-\sigma(j+2-k)} 
              \right),
$$
$$   B,~ \sigma \in{\bf C},~~~~~  -1 < \Re \sigma <1,~~~~~c_{jk}\in {\bf
C}$$ 
$B$ and $\sigma$ are integration constants. The $c_{jk}$'s are
	     certain  rational functions of $B$ and $\sigma$. For
	     symmetry reasons, we can restrict to the case 
$$
0\leq \Re\sigma <1.
$$ 
It
	     follows that: 
$$ 
\Omega_1^{(a)}(s) 
= b s^{-{\sigma\over 2}}\left( 1-{b^2
\over (1-\sigma)^2} s^{1-\sigma} 
 +{a^2 \over (1+\sigma)^2} s^{1+\sigma} +...\right)
 + a s^{\sigma \over 2} \left( 1+{b^2
\over (1-\sigma)^2} s^{1-\sigma} -{a^2 \over (1+\sigma)^2} s^{1+\sigma} + ...
\right),
$$
$$ 
\Omega_3^{(a)}(s) 
=i\left[ b s^{-{\sigma\over 2}}\left( 1-{b^2
\over (1-\sigma)^2} s^{1-\sigma}  +{a^2 \over (1+\sigma)^2} s^{1+\sigma} +
 ...\right) 
 - a s^{\sigma \over 2} \left( 1+{b^2
\over (1-\sigma)^2} s^{1-\sigma} -{a^2 \over (1+\sigma)^2} s^{1+\sigma} + ...
\right) \right],
$$
$$
 \Omega_2^{(a)}(s) 
= {i \sigma \over 2} + i{ b^2 \over 1 - \sigma}
 s^{1-\sigma} -i{a^2\over 1+\sigma}s^{1+\sigma}+...,
$$ 
where 
$$ 
b:= {R \over 2 B},~~~~~a:={R B \over 2},~~~~~~~\mu^2-{\sigma^2\over 4}= -R^2+o(1)
$$
 The dots are higher order
corrections, and the ordering depends on the  specific value of
$\sigma$. If we substitute into (\ref{bikkuri}) and we keep the
 dominant term, we obtain: 
$$
 y(s)\sim {4b^2\over (2\mu -\sigma)^2} s^{1-\sigma},~
~~~~~~~~~~~~~~~~~~~~~~~ \hbox{ for } 0<\Re
 \sigma<1
\hbox{ and for } \sigma=0,
$$
$$
 y(s)\sim s\Bigl[
\sin^2(\nu  \ln s + C(R,B)) +O(s) 
\Bigr],~~~~~~~~ \hbox{ for } \sigma=2i\nu,~~~\nu\in{\bf R}\backslash\{0\}.
$$
Thus, we have obtained the {\it
  small-power-type behaviors} {\bf 1)},   the Taylor series 
 degeneration of the {\it log-type
  behaviors} {\bf 4)} and the {\it sine-type oscillatory behavior}
 {\bf 2)}.

\vskip 0.2 cm 
{\bf ii)} Next, we consider the case corresponding to  $\sigma=1$. Let $\omega= \ln{\theta\over 4}+\gamma$, where
 $\gamma$ is the Euler's constant. In \cite{Tracy} we find the solution: 
$$
  \hat{y}(\theta)= -\theta\omega -{\theta^5\over
  128}(8\omega^3-8\omega^2+4\omega-1)+O(\theta^9\omega^5)
=
$$
$$
=
 -\theta \left(\ln {\theta\over 4} +\gamma\right)+ O(\theta^5 \ln^3
  \theta ),~~~~\theta>0,~~\theta\to 0
$$
   In the
variable $s$: 
$$
  \hat{y}(s)= i {R\over 2}  s^{1\over 2} \left(  \ln s + C
  \right)\Bigl[1+
O(s^2\ln^2s)\Bigr], ~~~~
C:=2\ln {R\over 4i}+2\gamma 
$$
 Now we compute:  
$$\Omega_2^{(a)}=  i s {1\over \hat{y}} {d \hat{y} \over ds} =
\left[ {i\over 2} 
+{i \over \ln s +C} + O(s^2 \ln s)\right]\Bigl(1+O(s^2\ln^2 s)\Bigr)
$$
$$
  \Omega_1^{(a)}(s)={R\over 2} (\hat{y}^{-1}+\hat{y})=\left[ -{i\over s^{1\over 2}(\ln s+C)} +
i{R^2\over 4} s^{1\over 2} (\ln s +C) \right] (1 +O(s^2 \ln^2s))
$$
$$
= -{i \over s^{1\over 2}(\ln s +C)} 
                                  +O(s^{1\over 2} \ln s)
$$
$$
 \Omega_3^{(a)}(s)= i{R\over 2} (y^{-1}-\hat{y}) = \left[ {1\over s^{1\over 2}
(\ln s+C)} + 
{R^2\over 4} s^{1\over 2} (\ln s +C) \right] (1 +O(s^2 \ln^2s)) 
$$
$$
= 
{1\over s^{1\over 2} (\ln s+C)} 
  +O(s^{1\over 2} \ln s)
$$
Observe that $\mu^2-{1\over 4}=-R^2+o(1)$.  Now,  
if we substitute into (\ref{bikkuri}) and we keep the
 dominant term, we obtain: 
$$
 y(s)\sim -{4\over (2\mu-1)^2} {1\over (\ln s +C)^2} 
$$
Thus, we have obtained a {\it log-type behavior} {\bf 4)}.

\vskip 0.2 cm 
{\bf iii)} The last case to be studied is:  
$$ 
  \sigma=1+2i \nu,~~~~\nu\in{\bf R}, ~~\nu\neq 0
$$
 In \cite{Tracy} we find the solution: 
$$
  y(\theta)=-{1\over 2\nu} ~\theta~\sin \left[2\nu \ln\left({\theta\over
  4}\right) + 2\varphi(\nu)\right]+O(\theta^3), ~~~\theta\to
  0,~~~~\varphi(\nu)=\arg\Gamma(i\nu). 
$$
 In variable $s$: 
$$
y(s)= i {R\over 2\nu} ~s^{1\over 2}~ \sin(\nu \ln s +D) +O(s^{3\over 2}) 
,~~~s\to 0,~~~~~~~D= 2 \nu\ln\left(-i{R\over 4}\right)+ 2\varphi(\nu)
$$
Therefore: 
$$
 \Omega_2^{(a)}= i s {1\over \hat{y}} {d \hat{y} \over ds} 
= {i \over \sin(\nu \ln s +D)
 +O(s)} \left[{1\over 2} \sin(\nu \ln s +D)+\nu \cos(\nu \ln s+D) +O(s)
 \right] =
$$
$$
={i\over 2} +i\nu{\cos(\nu \ln s+D)+O(s)\over \sin(\ln s + D)+O(s)},
$$

$$
\Omega_1^{(a)}= {R\over 2}(\hat{y}^{-1}+\hat{y})= -{1\over s^{1\over 2}}{i \nu\over 
\sin(\nu \ln s +D) +O(s)} +i{R^2\over 4 \nu} s^{1\over 2} (\sin(\nu \ln s +D)
+
O(s)),
$$
$$
\Omega_3^{(a)} = i{R\over 2}(\hat{y}^{-1}-\hat{y})={1\over s^{1\over
    2}} 
{\nu \over 
\sin(\nu \ln s +D) +O(s)} + {R^2\over 4 \nu} s^{1\over 2} (\sin(\nu \ln s + D)
+ O(s)). 
$$

\vskip 0.2 cm
\noindent
{\bf Remark:} 
 There is a sequence of poles (accumulating at $s=0$)  corresponding to the
roots of $\sin(\nu \ln s +D) +O(s)=0$. 
We stress that it is not possible to write $ {1 \over 
\sin(\nu \ln s +D) +O(s)}$ as  ${1 \over 
\sin(\nu \ln s +D) }(1 +O(s))$, because when we collect $\sin(\nu \ln s +D)$ in
the denominator we divide $O(s)$ by  $\sin(\nu \ln s +D)$ itself, so we
introduce poles (the roots of   $\sin(\nu
\ln s +D)=0$) in the $O(s)$ terms! 

\vskip 0.2 cm 
The computation of $y(s)$ is
 slightly more complicated than before. We have:  
$$
\Omega_3^{(a)}=i\Omega_1^{(a)}+O(s^{1\over 2}),
$$

$$
 A(s)\sim\left[
 {\Omega_1^{(a)}\Omega_2^{(a)}+i\mu\Omega_1^{(a)}\over 
\bigl(\Omega_2^{(a)}\bigr)^2+\mu^2} +{O(\sqrt{s}) \over R^2} 
\right]^2= \left[
{\Omega_1^{(a)}\over \Omega_2^{(a)} -i\mu} +O(\sqrt{s})
\right]^2=$$
$$
={1+O(\sqrt{s})\over 
s\left(
{2\mu-1\over 2\nu} \sin(\nu\ln s +D)-\cos(\nu \ln s+D) +O(s)
\right)^2
}
$$
Note that  in the computation we have used only the leading term of $\Omega_1^{(a)}$. The higher order terms $O(\sqrt{s})$ have been divided by $\bigl(\Omega_1^{(a)}\bigr)^2+\bigl(\Omega_3^{(a)}\bigr)^2=R^2$, while the leading part $\Omega_1^{(a)}\Omega_2^{(a)}+i\mu\Omega_1^{(a)}$ has been divided by $\bigl(\Omega_2^{(a)}\bigr)^2+\mu^2$, with the assumption that $\bigl(\Omega_1^{(a)}\bigr)^2+\bigl(\Omega_3^{(a)}\bigr)^2+\bigl(\Omega_2^{(a)}\bigr)^2\cong -\mu^2$.  
If we substitute the above result into (\ref{bikkuri}) and we keep the
 dominant term, we obtain: 

$$
y(s)\sim {1\over 1-\left({2\mu -1 \over 2\nu} \sin(\nu\ln s +D) -\cos(\nu
  \ln s+D)\right)^2+O(\sqrt{s})}
=
$$

$$
=
{1\over 1 -{4\nu^2+(2\mu-1)^2\over 4\nu^2} \sin^2(\nu \ln s +C)+O(\sqrt{s})},~~~~~C:=D-{i\over 2} \ln{2\mu-1-2i\nu\over 2\mu-1+2i\nu}. 
$$
Thus, we have found an {\it inverse sine-type oscillatory behavior}
{\bf 3)}.

\section{On the error in the approximation of PVI with PIII}  
We may estimate the error in the approximation of $\Omega_j$ with
$\Omega_j^{(a)}$.   Let's 
estinate the error in $\Omega_3$. We put
$$
 \Omega_3: =\Omega_3^{(a)}+\delta \Omega_3
$$ 
and we substitute into the left
 hand side of the third equation of 
 (\ref{system}). In the right hand side we substitute $\Omega_1^{(a)}$ and
 $\Omega_2^{(a)}$:
$$ 
  {d \Omega_3^{(a)} \over ds} + {d \over ds} (\delta \Omega_3)=
     -{1\over s} ( 1 +s +...) ~\Omega_1^{(a)}\Omega_2^{(a)}
$$
Recalling that $ {d \Omega_3^{(a)} \over ds}=-{1\over
s}~\Omega_1^{(a)}\Omega_2^{(a)}$ we get 
$$
 {d \over ds} (\delta \Omega_3)=-\Omega_1^{(a)}\Omega_2^{(a)}+...\sim -i{\sigma
\over 2}  (bs^{-{\sigma\over 2}}+ a s^{\sigma\over 2})+...
$$
 which implies 
$$
  \delta \Omega_3 \cong O(s^{1-{\sigma\over 2}}+s^{1+{\sigma\over 2}})
$$
We write $\Omega_1:= \Omega_1^{(a)} + \delta \Omega_1$ and we substitute in the
 left hand side of the first equation of (\ref{system}), while in the right
 hand side we substitute the  $\Omega_2^{(a)}$, $ \Omega_3^{(a)}$. The same
 procedure yields:
$$
  \delta \Omega_1 \cong O(s^{1-{\sigma\over 2}}+s^{1+{\sigma\over 2}})
$$
In fact,  the terms of order $s^{1\pm {\sigma\over 2}}$ are
missing in the approximated solutions $\Omega_1^{(a)}$ and
$\Omega_3^{(a)}$,
 but they appear in the true formal expansion of the $\Omega_1$ and
 $\Omega_2$, which is computed in subsection \ref{small}. As for $\Omega_2$, we proceed as above
making use of the second equation of (\ref{system}), which becomes:
$$
 {d \delta \Omega_2 \over ds}= s \Omega_1^{(a)}\Omega_3^{(a)}
= i b^2 s^{1-\sigma} -i a^2 s^{1+\sigma} +...
$$
Thus 
$$
\delta \Omega_2= i{b^2\over {2-\sigma} } s^{2-\sigma} - i{a^2 \over
2+\sigma} s^{2+\sigma}+\hbox{ higher orders}
$$
}

\subsection{Expansion with respect to a Small Parameter}
\label{small} 
The error in the asymptotic reduction is more precisely evaluated if
we write the true formal expansion of the $\Omega_j$'s. In order to do this,   a small parameter expansion can be used, as it is done in \cite{guz2001}.  
 Let $s:=\epsilon z$ 
where $\epsilon$ is the small
parameter. The system (\ref{system}) becomes: 
\be
{d \Omega_1\over dz} = {1 \over z} ~ \Omega_2 \Omega_3,~~~
{d \Omega_2\over dz} = {\epsilon \over 1-\epsilon z} ~ \Omega_1 \Omega_3,~~~
{d \Omega_3\over dz} = {1 \over z( \epsilon z -1)} ~ \Omega_1
\Omega_2.   
\ee
The coefficient of the new system are holomorphic for $\epsilon \in E:=
\{ \epsilon\in {\bf C}~|~|\epsilon| \leq \epsilon_0 \}$ and for $0<|z|<{1\over
|\epsilon_0|}$, in  
particular for $z\in D:=\{z\in {\bf C} |~R_1\leq |z| \leq R_2 \}$,
 where $R_1$ and $R_2$ are independent of $\epsilon$  and satisfy 
$0<R_1<R_2<{1\over \epsilon_0}$. The small parameter expansion is a  formal way to compute the
 expansions of the $\Omega_j$'s  for $s\to 0$. 
To our knowledge, the procedure does not give a  rigorous 
  justification of the uniform  
convergence of the $s$-expansions of the $\Omega_j$'s. 
 For $\epsilon \in E$ and $z\in D$ we can expand the fractions as
 follows: 
\be
{d \Omega_1\over dz} = {1 \over z} ~ \Omega_2 \Omega_3,~~~
{d \Omega_2\over dz} = \epsilon ~\sum_{n=0}^{\infty} z^n \epsilon^n 
 ~ \Omega_1 \Omega_3,~~~
{d \Omega_3\over dz} = -{1 \over z}~\sum_{n=0}^{\infty} z^n \epsilon^n
 ~ \Omega_1 \Omega_2  
\label{systpower}
\ee
and we look for a solution  expanded in powers of $\epsilon$:
\be
   \Omega_j(z,\epsilon)=\sum_{n=0}^{\infty}
   \Omega_j^{(n)}(z)~\epsilon^n~,~~~~
~~j=1,2,3.
\label{solution}
\ee
  We find the $\Omega_j^{(n)}$'s substituting (\ref{solution}) into
  (\ref{systpower}). At order $\epsilon^0$ we find
$$
\left\{ \matrix{ 
{\Omega_2^{(0)}}^{\prime}=0~~~~~~~~~~~~ \cr
\cr
{\Omega_1^{(0)}}^{\prime}={1\over z}~ \Omega_2^{(0)}\Omega_3^{(0)} 
\cr
\cr
 {\Omega_3^{(0)}}^{\prime}=-{1\over z}~ \Omega_2^{(0)}\Omega_1^{(0)} 
}\right.
$$
The prime denotes
the derivative w.r.t. $z$. Thus:
$$
\Omega_2^{(0)}={i \sigma
\over 2}
$$
 Then we solve the linear system for
$\Omega_1^{(0)}$ and $\Omega_3^{(0)}$ and find
$$ 
   \Omega_1^{(0)}= \tilde{b}~ z^{-{\sigma \over 2}}+\tilde{a} 
~ z^{\sigma\over 2}= (\tilde{b}\epsilon^{\sigma\over 2})
   ~s^{-{\sigma\over 2}} +(\tilde{a} \epsilon^{-{\sigma\over
   2}})s^{\sigma\over 2}
$$
$$
   \Omega_3^{(0)}= i\tilde{b}~z^{-{\sigma \over 2}}-i\tilde{a}
~ z^{\sigma\over 2}= 
i(\tilde{b}\epsilon^{\sigma\over 2})
   ~s^{-{\sigma\over 2}} -i(\tilde{a} \epsilon^{-{\sigma\over
   2}})s^{\sigma\over 2}
$$
where $\tilde{a}$ and $\tilde{b}$ are integration constants. We will
require that $b:=\tilde{b}\epsilon^{\sigma\over 2}$,  $a:=\tilde{a} \epsilon^{-{\sigma\over
   2}}$ are finite, when $\epsilon \to 0$. 
The higher orders are: 
$$
\left\{ 
\matrix{
\Omega_2^{(n)}(z)=\int^z d \zeta ~
\sum_{k=0}^{n-1}~\zeta^k\sum_{l=0}^{n-1-k}~\Omega_1^{(l)}(\zeta)
~\Omega_3^{(n-1-k-l)} (\zeta)\cr\cr
 {\Omega_1^{(n)}}^{\prime}={1\over z}~ \Omega_2^{(0)}\Omega_3^{(n)}+
A_1^{(n)}(z)
\cr
\cr 
 {\Omega_3^{(n)}}^{\prime}=-{1\over z}~
\Omega_2^{(0)}\Omega_1^{(n)}+A_3^{(n)}(z)
}
\right.
$$
where:  
$$
\left\{
\matrix{
A_1^{(n)}(z)= {1\over z} \sum_{k=1}^n ~\Omega_2^{(k)}(z)~\Omega_3^{(n-k)}(z)
\cr
\cr
A_3(z)= -{1\over z} \left[
                          \sum_{l=1}^n~\Omega_2^{(l)}(z)~ \Omega_1^{(n-l)}(z) +
\sum_{k=1}^n ~z^k\sum_{l=0}^{n-k} ~\Omega_1^{(l)}(z)~ \Omega_2^{(n-k-l)}(z)
                          \right] 
}
\right.
$$
The system for $\Omega_1^{(n)}$, $\Omega_3^{(n)}$ is closed and
non-homogeneous. By variation of parameters we find the particular
solution 
$$
\Omega_1^{(n)}(z)= {z^{\sigma/2} \over \sigma} ~\int^{z} d \zeta~
\zeta^{1-{\sigma\over 2}} R_1^{(n)}(\zeta)~-{z^{-{\sigma/2}} \over
\sigma} ~\int^z \zeta^{1+{\sigma\over 2}} R_1^{(n)}(\zeta)
$$
$$
\Omega_3^{(n)}(z)= {z \over i \sigma /2} \left(
{\Omega_1^{(n)}}(z)^{\prime}-A_1^{(n)}(z)   \right)
$$
where
$$ R_1^{(1)}(z)= {1\over z} A_1^{(n)}(z) + {i \sigma \over 2 z} 
A_3^{(n)}(z) + {A_1^{(n)}}(z)^{\prime}
$$
Thus: 

\be
\left\{\matrix{
\Omega_j(s)= s^{-{\sigma \over 2}} \sum_{k,~q=0}^{\infty} b_{kq}^{(j)} 
~s^{k +(1-\sigma)q} + s^{{\sigma \over 2}} \sum_{k,~q=0}^{\infty} a_{kq}^{(j)} 
~s^{k +(1+\sigma)q} & j=1,3\cr\cr
\Omega_2=  \sum_{k,~q=0}^{\infty} b_{kq}^{(2)} 
~s^{k +(1-\sigma)q} + \sum_{k,~q=0}^{\infty} a_{kq}^{(2)} 
~s^{k +(1+\sigma)q}&\cr
}
\right.
\label{smallserie}
\ee

\noindent
The coefficients $a_{kq}^{(j)}$ and $b_{kq}^{(j)}$ contain
$\epsilon$. 
In fact, they are functions of $a:=\tilde{a} 
\epsilon^{-{\sigma\over 2}}$, $b:=\tilde{b} \epsilon^{{\sigma\over
2}}$. The re-normalization after restoring $s$  is possible if 
$\sigma\neq 2n+1$, $n\in{\bf Z}$ (no $\ln
z$ terms in $\Omega_2$), and if the additive constant in the
integration of $\Omega_2^{(n)}$ is zero. If this is not the case, some 
 coefficients of the expansions for the $\Omega_j$'s diverge.

\vskip 0.2 cm 
 We can fix the range of $\sigma$ according to the condition that the first term $\Omega_2^{(0)}$ in $\Omega_2$ be the leading one. 
The 
approximation at order 0 for $\Omega_2$ is:
$$
  \Omega_2 \approx {i \sigma \over 2} \equiv \hbox{ constant}
$$
The approximation at order 1 contains powers $z^{1-\sigma}$,
$z^{1+\sigma}$. If we assume that the  approximation at order 0 in
$\epsilon$ is actually the limit of $\Omega_2$ as $s=\epsilon~z \to 0$,
than we need 
$$ 
  -1 < \Re \sigma <1
$$
The ordering of the expansion (\ref{smallserie}) is
 somehow conventional: namely, we could transfer some terms multiplied
 by $s^{{\sigma\over 2}}$ in the series multiplied by
 $s^{-{\sigma\over 2}}$, and conversely.  I report
 the first terms: 
$$
\Omega_1(s) = b s^{-{\sigma\over 2}} \left(1 -{b^2\over
(1-\sigma)^2}s^{1-\sigma}+{\sigma^2\over 4(1-\sigma)}s+ {a^2\over
(1+\sigma)^2}s^{1+\sigma} 
+...\right)$$
$$
+as^{\sigma\over 2}\left(1+{b^2\over
(1-\sigma)^2}s^{1-\sigma}+{\sigma^2\over 4(1+\sigma)}s-  {a^2\over
(1+\sigma)^2}s^{1+\sigma} 
+... \right)
$$
\vskip 0.2 cm 
$$
\Omega_3(s)=ibs^{-{\sigma\over 2}}\left(
           1-{b^2\over (1-\sigma)^2}s^{1-\sigma}+{\sigma(\sigma-2) \over
           4(1-\sigma)} s +{a^2\over (1+\sigma)^2}s^{1+\sigma} +...
\right)$$
$$
-ias^{\sigma\over 2} \left( 
           1+{b^2\over (1-\sigma)^2} s^{1-\sigma} +{\sigma (\sigma+2) \over
           4(1+\sigma)} s - {a^2 \over (1+\sigma)^2} s^{1+\sigma} +...
\right)
$$
\vskip 0.2 cm 
$$
\Omega_2(s)= i~{\sigma\over 2}+i {b^2\over 1-\sigma} s^{1-\sigma}
  -i {a^2\over 1+\sigma} s^{1+\sigma} +...
$$
Note that the dots do not mean higher order terms. There may be 
 terms bigger than
 those written above (which are computed through the expansion in the small
 parameter up to order $\epsilon$) depending on the particular value of $\Re \sigma $
 in $(-1,1)$. 
Finally, we note that we  can always assume:
$$
   0\leq \Re \sigma <1, 
$$ 
because that would not affect the expansion of the solutions but for
the change of two signs.

\vskip 0.2 cm
It is worth observing that if we pretend that the
solutions (\ref{smallserie}) are still valid for  $ \Re \sigma =1$  and if we
extract the terms where $s$ has exponent with negative or vanishing 
real part, we have: 
$$
 \Omega_1 = b s^{-{1\over 2} -i \nu} ~\sum_{q=0}^{\infty} 
(-1)^q \left[{b^2 \over
 (1-\sigma)^2}\right]^q ~s^{(1-\sigma)q}|_{\sigma=1+2 i \nu}+...
$$
  $$
 \Omega_3 =i b s^{-{1\over 2}- i \nu} ~\sum_{q=0}^{\infty} (-1)^q \left[{b^2 \over
 (1-\sigma)^2}\right]^q ~s^{(1-\sigma)q}|_{\sigma=1+2 i \nu}+...
$$ 
$$
 \Omega_2= {i \sigma \over 2} + i {b^2 \over 1-\sigma} ~s^{-2i\nu}
\sum_{q=0}^{\infty} ~(-1)^q\left({b^2\over (1-\sigma)^2}\right)^q
~s^{(1-\sigma)q}|_{\sigma=1+2 i \nu} +...
$$
For $|\Re \sigma|<1$ an $s$ small, we sum the series:  
$$ 
 \sum_{q=0}^{\infty} ~(-1)^q \left[\left({b\over 1-\sigma}\right)^2
 s^{1-\sigma}   \right]^q= {1\over 1 + \left({b\over 1-\sigma}\right)^2
 s^{1-\sigma} } 
$$
 Then, we analytically extend the result at $\sigma=1+2i\nu$. In this way 
$$
 \Omega_1= -{i\nu \over s^{1\over 2} \sin(\nu \ln s +C)} +...$$
$$
 \Omega_3= {\nu \over s^{1\over 2} \sin(\nu \ln s +C)} +...$$
$$
 \Omega_2= {i\over 2} +i \nu \cot(\nu \ln s +C) +...
$$
 where $C= -i \ln (2\nu/b)$. The result is similar to that obtained
 for the $\Omega_j^{(a)}$'s, but the $O(s)$-terms in the denominator do not appear.

\vskip 0.3 cm 
\noindent 
{\bf NOTE:} {\small 
  Consider the system (\ref{system}) and expand the fractions as $s\to
 0$. We find
\be
{d \Omega_1\over ds} = {1\over s} ~ \Omega_2 \Omega_3,~~~
{d \Omega_2\over ds} = \sum_{n=0}^{\infty} s^n ~ \Omega_1 \Omega_3,~~~
{d \Omega_3\over ds} = - {1 \over s}\sum_{n=0}^{\infty} s^n
                       ~ \Omega_1 \Omega_2  
\ee
  We can look for a formal  solution: 
$$
\Omega_j(s)= s^{-{\sigma \over 2}} \sum_{k,~q=0}^{\infty} b_{kq}^{(j)} 
~s^{k +(1-\sigma)q} + s^{{\sigma \over 2}} \sum_{k,~q=0}^{\infty} a_{kq}^{(j)} 
~s^{k +(1+\sigma)q}~,~~~~~~~j=1,3
$$
$$
\Omega_2=  \sum_{k,~q=0}^{\infty} b_{kq}^{(2)} 
~s^{k +(1-\sigma)q} + \sum_{k,~q=0}^{\infty} a_{kq}^{(2)} 
~s^{k +(1+\sigma)q}
$$
 Plugging the series into the equation we find solvable relations
 between the coefficients and we can determine them. For example, the
 first relations  give 
$$
  \Omega_2 = {i \sigma \over 2} + \left({i [b_{00}^{(1)}]^2 \over 1 -\sigma}
  s^{1-\sigma}+...\right) - \left({i [a_{00}^{(1)}]^2 \over 1 +\sigma}
  s^{1+\sigma}+...\right)$$
$$
\Omega_1= (b_{00}^{(1)}s^{-{\sigma\over 2}}+...)+ (a_{00}^{(1)}
s^{\sigma \over 2} +...)$$
$$
\Omega_3= (ib_{00}^{(1)}s^{-{\sigma\over 2}}+...)+ (-ia_{00}^{(1)}
s^{\sigma \over 2} +...)$$
All the  coefficients determined by successive relations  
are functions of $\sigma$, $b_{00}^{(1)}$, $
a_{00}^{(1)}$. These are the three parameters on which the solution of
 (\ref{system}) must depend.  We can identify $b_{00}^{(1)}$ with $b$ and
$a_{00}^{(1)}$ with $a$. 
}


\vskip 0.3 cm
The case  $\sigma=1$ in the small parameter formalism is more
complicated. 
If we perform the small parameter
expansions as before, we find the same $\Omega_j^{(0)}$ than
before. But due to the exponent $z^{-1/2}$ the integration for
$\Omega_2^{(1)}$ gives 
$$
 \Omega_2^{(1)} = -{i\over 2} \tilde{a}^2 z^2 + i \tilde{b}^2\ln(z)
$$
 In this way, we find for $\Omega_1$ and $\Omega_3$ an
 expansion in power of $\epsilon$ with  coefficients which are
 polynomials in $\ln(z)$; also the powers $z^{-1/2}$, $z^{1/2}$,... ,
 $z^{n/2}$, $n>0$ appear in the coefficients. $\Omega_2$ is an 
 expansion in power of $\epsilon$ with coefficients which are
 polynomials in $\ln(z)$ and $z$. It is not obvious how to recombine $z$ and
$\epsilon$ when logarithms appear. We can put $\epsilon=1$. Anyway, we
see that the first correction to the constant ${i \over 2}$ in
$\Omega_2$ is  $\ln(s)$, which is {\it not a correction} to the
constant when $s\to 0$, because it diverges. 

\vskip 0.2 cm

 We can try an expansion which already contains logarithms of
 $\epsilon$: 
\be
   \Omega_j(z,\epsilon)= \sum_{k=-1}^{+\infty}  \sum_{n=0}^{+\infty}
    \Omega_{j,n}^{(k)}(z) ~{\epsilon^{2 k +1 \over 2} \over (\ln
 \epsilon)^n}
,~~~~j=1,3
\label{chazy1}
\ee
\be
 \Omega_2=  \sum_{k=0}^{+\infty}  \sum_{n=0}^{+\infty} 
\Omega_{2,n}^{(k)}(z) ~{\epsilon^{k} \over (\ln
 \epsilon)^n}
\label{chazy2}
\ee
Then we substitute in (\ref{systpower}) and we equate powers of
$\epsilon$ and $\ln \epsilon$. The
requirement that we could re-compose the powers of $\ln \epsilon$ and
$\ln z$ appearing in the expansion in the form $\ln (z \epsilon)$
imposes very strong relations on the integration constants. The result
 which we obtain, when we solve the equations for the   
coefficients $\Omega_{j,n}^{(k)}$ equating powers up to ${1\over( \ln
\epsilon)^n}$ and $\epsilon^{-1/2}$, is : 
$$
 \Omega_1 = {i \over  s^{1\over2}(\ln(s) + C)} +O\left({1\over (\ln
 \epsilon)^n}\right)
               + O(\epsilon^{1\over 2}) ,
$$
$$
 \Omega_3 = {-1 \over  s^{1\over2}(\ln(s) + C)} +O\left({1\over (\ln
 \epsilon)^n}\right)
               + O(\epsilon^{1\over 2}) ,
$$
$$ 
\Omega_2 = {i \over 2} + {i \over \ln s + C } + O\left({1\over (\ln
 \epsilon)^n}\right)
               + O(\epsilon^{1\over 2}).  
$$ 
  These are the log-type behaviors.

\section{Conclusions}
 The  asymptotic reduction of PVI to PIII  produces the
 correct leading term of the critical
 behavior of branches of PVI transcendents, classified in \cite{guz2010}, starting from the asymptotic
 behaviors of PIII transcendents, computed in \cite{Tracy}. We have evaluated the error of the asymptotic reduction, showing that it does not affect the leading term of the PVI transcendents obtained from asymptotic behaviors of PIII transcendents. 
 

\end{document}